\documentclass[11pt]{article}
\usepackage[T1]{fontenc}
\usepackage{lmodern}
\usepackage[margin=1.08in]{geometry}
\usepackage{amsmath,amssymb,amsthm,mathtools}
\numberwithin{equation}{section}
\usepackage{microtype}
\usepackage{enumitem}
\usepackage{xcolor}
\usepackage[pdftitle={Proper sample compression schemes for VC classes},pdfauthor={Artem Chernikov},colorlinks=true,linkcolor=blue!45!black,citecolor=blue!45!black,urlcolor=blue!45!black]{hyperref}
\setlist{itemsep=2pt,topsep=3pt}
\newtheorem{theorem}{Theorem}[section]
\newtheorem{fact}[theorem]{Fact}
\newtheorem{lemma}[theorem]{Lemma}

\theoremstyle{definition}
\newtheorem{remark}[theorem]{Remark}
\newtheorem{problem}[theorem]{Problem}
\newtheorem{definition}[theorem]{Definition}
\newcommand{\constof}[1]{C_{\ref*{#1}}}
\DeclareMathOperator{\VC}{VC}
\DeclareMathOperator{\dom}{dom}
\newcommand{\C}{\mathcal C}
\newcommand{\Z}{\mathcal Z}
\newcommand{\E}{\mathcal E}
\newcommand{\bits}{\{0,1\}}
\title{Proper sample compression schemes for VC classes}
\author{Artem Chernikov}
\date{}
\begin{document}
\maketitle
\vspace{-1.3em}
\begin{abstract}
Every nonempty binary concept class of VC dimension at most $d$ admits a labeled sample compression scheme of total size $2^{O(d)}$ which is additionally proper, i.e., its reconstruction map always returns a member of the original class.
\end{abstract}

\section{Introduction}
Sample compression, introduced by Littlestone and Warmuth \cite{LW}, asks how much of a labeled sample one must retain to recover all its labels. A compression scheme specifies a small subsample, possibly together with a bounded amount of side information, from which a fixed reconstruction map produces a hypothesis consistent with the original sample. Under the usual measurability assumptions, such schemes give generalization bounds in the PAC framework \cite{Valiant,FW}. Motivated by the VC characterization of binary PAC learnability \cite{VC,BEHW}, Littlestone and Warmuth asked whether every class of VC dimension $d$ admits a compression scheme of size bounded solely in terms of $d$.

Early boosting schemes retained $O(d\log n)$ examples from a sample of size $n$ \cite{Freund,HelmboldWarmuth,FW}. By finitizing a model-theoretic result of Chernikov and Simon \cite{CS}, Livni and Simon \cite{LS} gave a bound independent of $n$ assuming that certain classes obtained from the original incidence relation by Boolean operations and projection have finite VC dimension. For a class of $M$ concepts and fixed $d$, Moran, Shpilka, Wigderson, and Yehudayoff \cite{MSWY} retained $O_d(\log\log M)$ examples. Moran and Yehudayoff \cite{MY} obtained the first bound depending only on $d$ for arbitrary VC classes, showing that total size $2^{O(d)}$ suffices. Cai et al.~\cite{Cai} recently provided a nearly quadratic bound in $d$ for ordinary compression. The linear bound conjectured by Floyd and Warmuth \cite{FW} remains open.

Here we consider \emph{proper} compression, where the reconstructed hypothesis belongs to the given class. Several special cases were known on finite domains. Proper labeled schemes of size $d$ are known for extremal classes \cite{MW} and for topes of complexes of oriented matroids \cite{CKP}. Maximum classes admit schemes of size $d$ \cite{FW}, even with the labels omitted \cite{CCMW}. Balls of arbitrary radius in trees, cycles, interval graphs, trees of cycles, and cube-free median graphs admit constant-size proper labeled schemes \cite{CCMRV}. In any finite graph of treewidth or clique-width at most $t$, the class of balls admits a proper scheme of size $O(t\log t)$ \cite{BHMR}. Below we mostly follow the notation of \cite[Section 1.2]{MY}.

\begin{definition}\label{def:scheme}
Let $\varnothing\ne\C\subseteq\bits^X$. For $B\subseteq X$, write $\C|_B\coloneqq \{c|_B:c\in\C\}$.
\begin{enumerate}[label=\textup{(\roman*)},leftmargin=*]
\item A finite set $B\subseteq X$ is \emph{shattered} if $\C|_B=\bits^B$. The VC dimension $\VC(\C)$ is the largest size of a finite shattered set, or $\infty$ if these sizes are unbounded.
\item A \emph{sample} is a pair $(B,\ell)$ with $B\subseteq X$ finite and $\ell:B\to\bits$. It is \emph{realizable} if $\ell\in\C|_B$. Write $L_\C(\infty)$ for the set of realizable samples and $L_\C(k)$ for those with $|B|\le k$.
\item A \emph{labeled sample compression scheme} with kernel size $k$ and finite nonempty side-information set $I$ consists of a compression map $\kappa$ and a reconstruction map $\rho$
\[
 \kappa:L_\C(\infty)\longrightarrow L_\C(k)\times I,
 \qquad \rho:L_\C(k)\times I\longrightarrow\bits^X,
\]
so that for each $(B,\ell)\in L_\C(\infty)$ there is a subset $T\subseteq B$ of size at most $k$ and a tag $i\in I$ such that
\[
 \kappa(B,\ell)=((T,\ell|_T),i),\qquad \rho(\kappa(B,\ell))|_B=\ell.
\]
\item The scheme is \emph{proper} if the range of $\rho$ is contained in $\C$.
\item The \emph{total size} of the scheme is $k+\log_2|I|$.
\end{enumerate}
\end{definition}

Answering a question of Ben-David et al.~\cite[Section 3]{BDHMSY}, we prove the following.

\begin{theorem}\label{thm:main}
There is an absolute constant $K>0$ such that for every integer $d\ge 0$ and every set $X$, every nonempty binary concept class $\C\subseteq\bits^X$ of VC dimension at most $d$ admits a proper labeled sample compression scheme of total size at most $2^{Kd}$.
\end{theorem}

This also gives a monotone compression scheme for $\C$ of size $2^{O(d)}$ in the sense of \cite[Section~3]{BDHMSY} (using the equivalence observed there).

\paragraph{Proof outline.}
We first describe the proof for finite $X$. Section~\ref{sec:arbitrary} extends the construction to arbitrary $X$ with the same bounds by a compactness argument. Moran and Yehudayoff \cite{MY} reconstruct by majority vote. In our proof a majority vote gives an intermediate partial labeling, and the reconstruction map returns a member of $\C$. For each nonempty subclass $\Z\subseteq\C$, we construct two inclusion-preserving maps $R_\Z:2^X\to 2^\Z$ and $P_\Z:2^X\to 2^X$ (Definition~\ref{def:descriptions}). For $T\subseteq X$, the sets $R_\Z(T)$ and $P_\Z(T)$ consist of the candidate concepts and test points available from $T$, respectively. By Lemma~\ref{lem:availability}, each such concept or point is already available from some $T_0\subseteq T$ of size bounded in terms of $d$, and both $|R_\Z(T)|$ and $|P_\Z(T)|$ are bounded in terms of $|T|$ and $d$.

The main step is Lemma~\ref{lem:dimension-descent}. It gives the following alternative for $a\in\Z$ and $B\subseteq X$ with $r\coloneqq\VC(\Z|_B)>0$. Either some $c\in R_\Z(B)$ agrees with $a$ on $B$, or there is a set $F\subseteq P_\Z(B)$ of size bounded in terms of $d$ such that $\{c\in\Z:c|_F=a|_F\}$ has VC dimension strictly less than $r$ on $B$. If no such $F$ exists, Lemmas~\ref{lem:maj} and~\ref{lem:support} give a partial labeling $h$ agreeing with $a$ on $B$ and at every point where Lemma~\ref{lem:correct} tests $h$ for disagreement. Lemma~\ref{lem:correct} then gives the desired $c$. At dimension zero, $\min\Z\in R_\Z(B)$ already agrees with $a$ on $B$.

Starting with $\Z=\C$, at most $d$ such restrictions, each reducing the VC dimension on $B$, suffice to obtain a candidate agreeing with the sample. By Lemma~\ref{lem:availability}, retaining boundedly many points for each test point used in the restrictions and for the final candidate gives a set $T\subseteq B$ of size $2^{O(d)}$ (Lemma~\ref{lem:retained-set}). The decoder enumerates all admissible restrictions and their candidate concepts using $T$ alone (Definition~\ref{def:candidate-tree}). By Lemma~\ref{lem:candidate-count}, the position of the desired concept in this list is specified by $2^{O(d)}$ bits.

An important ingredient in Lemma~\ref{lem:correct} is a polynomial strengthening of Matou\v{s}ek's $(p,q)$-theorem for VC classes (Lemma~\ref{lem:pq}). For infinite $X$, we also use the infinitary $(p,q)$-theorem of Jung and P\'alv\"olgyi \cite[Theorem~4]{JP}. Matou\v{s}ek's qualitative theorem underlies the aforementioned compression argument in theories without the independence property (NIP) \cite{CS}, and plays an important role in further model-theoretic developments \cite{EK,KaplanPQ,BKS,BKR}.

\section{A polynomial \texorpdfstring{$(p,q)$}{(p,q)}-theorem for VC families}\label{sec:estimates}
\begin{definition}
Until Section~\ref{sec:arbitrary}, fix a finite set $X$, an integer $d\ge1$, and a nonempty class $\C\subseteq\bits^X$ with $\VC(\C)\le d$.
\end{definition}

Write $[j]\coloneqq\{1,\ldots,j\}$. For finite sets $I$ and $Y$, an \emph{indexed family} $\E\coloneqq(E_i)_{i\in I}$ of subsets of $Y$ assigns a set $E_i\subseteq Y$ to each $i\in I$, allowing repetitions. Its dual is the indexed family $\E^*\coloneqq(R_y)_{y\in Y}$ of subsets of $I$, where $R_y\coloneqq\{i\in I:y\in E_i\}$. For an indexed family, the VC dimension and shatter function refer to those of the underlying family of distinct sets. Let $\pi_{\E}$ be the shatter function of $\E$, i.e., $\pi_{\E}(t)$ is the maximum of $|\{E_i\cap B:i\in I\}|$ over $B\subseteq Y$ with $|B|\le t$. Let $\pi^*_{\E}\coloneqq\pi_{\E^*}$ denote the dual shatter function of $\E$. The following is standard.

\begin{lemma}\label{lem:vc-bounds}
Let $\E\coloneqq(E_i)_{i\in I}$ be an indexed family of subsets of $Y$, with $I,Y$ finite and nonempty and $\VC(\E)\le s$ for an integer $s\ge1$.
\begin{enumerate}[label=\textup{(\roman*)},leftmargin=*,nosep]
\item For $H\subseteq A\subseteq Y$, put $D_i\coloneqq(E_i\cap A)\mathbin{\triangle}H$ for $i\in I$. The indexed family $\mathcal D\coloneqq(D_i)_{i\in I}$ shatters precisely the subsets of $A$ shattered by $\E$. In particular, $\VC(\mathcal D)\le\VC(\E)$.
\item Put $D\coloneqq2^{s+1}-1$. Then $\VC(\E^*)\le D$, $\pi_{\E}(t)\le(t+1)^s$ and $\pi^*_{\E}(t)\le(t+1)^D$ for $t\ge1$.
\item For an integer $m\ge1$, let $\mathcal U\coloneqq(U_\alpha)_{\alpha\in I^m}$ be the indexed family of all $m$-fold unions, where $U_\alpha\coloneqq\bigcup_{k=1}^m E_{\alpha_k}$ for $\alpha=(\alpha_1,\ldots,\alpha_m)\in I^m$. If $\VC(\E^*)\le s^*$ and $s^*\ge1$, then $\VC(\mathcal U)=O(ms\log(2m))$ and $\VC(\mathcal U^*)=O(s^*\log(2m))$.
\end{enumerate}
\end{lemma}
\begin{proof}
For $B\subseteq A$, we have $D_i\cap B=(E_i\cap B)\mathbin{\triangle}(H\cap B)$. Taking the symmetric difference with $H\cap B$ permutes the subsets of $B$, proving (i). Part (ii) follows from Assouad's dual VC bound \cite[Proposition~2.12]{Assouad} and the Sauer--Shelah lemma \cite{Sauer,Shelah} applied to $\E$ and $\E^*$.

For (iii), taking unions of traces gives $\pi_{\mathcal U}(t)\le\pi_{\E}(t)^m$. Fix $J\subseteq I^m$ of size at most $t$ and put $W\coloneqq\{\alpha_k:\alpha\in J,\ k\in[m]\}$. For $y\in Y$, its trace $A_y\coloneqq\{i\in W:y\in E_i\}$ determines $\{\alpha\in J:y\in U_\alpha\}=\{\alpha\in J:\alpha_k\in A_y\text{ for some }k\in[m]\}$. Since $|W|\le mt$, we have $\pi^*_{\mathcal U}(t)\le\pi^*_{\E}(mt)$.

If $\mathcal U$ shatters $t\ge ms$ points, Sauer--Shelah gives $2^t\le(et/s)^{ms}$. If $\mathcal U^*$ shatters $t\ge s^*$ indices, it gives $2^t\le(emt/s^*)^{s^*}$. Put $z\coloneqq t/(ms)$ or $z\coloneqq t/s^*$, respectively. Then $2^z\le emz$, so $z=O(\log(2m))$ since $\log_2 z\le z/2$ for $z\ge4$. This proves (iii).
\end{proof}

For a finite nonempty set $S$, let $\Delta(S)$ denote the set of probability distributions on $S$.

\begin{fact}[von Neumann's minimax theorem \cite{vonNeumann}]\label{fact:minimax}
Let $M$ be a real matrix with finite nonempty row and column index sets $I$ and $Y$. Then
\[
 \min_{\mu\in\Delta(I)}\max_{y\in Y}\sum_{i\in I}\mu(i)M_{iy}
   =\max_{\nu\in\Delta(Y)}\min_{i\in I}\sum_{y\in Y}\nu(y)M_{iy}.
\]
The common value is called the \emph{value} of $M$.
\end{fact}

Moran and Yehudayoff \cite{MY} use the following consequence of the $\varepsilon$-approximation theorem of Vapnik and Chervonenkis \cite{VC}. Fact~\ref{fact:vc-minimax} uses the improved bound in \cite{Talagrand,LLS}. We write $C_\nu$ for a constant supplied by statement $\nu$. When it depends on $d$, we write $C_\nu(d)$ and suppress the argument once $d$ is fixed.

\begin{fact}\label{fact:vc-minimax}
Let $\E\coloneqq(E_i)_{i\in I}$ be an indexed family of subsets of $Y$, with $I,Y$ finite and nonempty. Write $R_y\coloneqq\{i\in I:y\in E_i\}$ for $y\in Y$. The incidence matrix $M$ of $\E$ has rows $i\in I$, columns $y\in Y$, and entries $M_{iy}\coloneqq\mathbf{1}_{\{y\in E_i\}}$. Let $v$ be its value, and put $s\coloneqq\VC(\E)$ and $s^*\coloneqq\VC(\E^*)$. There is an absolute constant $\constof{fact:vc-minimax}>0$ such that \textup{(i)} and \textup{(ii)} hold for every $0<\varepsilon<1$.
\begin{enumerate}[label=\textup{(\roman*)}]
\item There are columns $(y_j)_{j\in[m]}\in Y^m$, with $1\le m\le \constof{fact:vc-minimax}(s+1)/\varepsilon^2$, such that $|\{j\in[m]:y_j\in E_i\}|/m\ge v-\varepsilon$ for every $i\in I$.
\item There are rows $(i_j)_{j\in[n]}\in I^n$, with $1\le n\le \constof{fact:vc-minimax}(s^*+1)/\varepsilon^2$, such that $|\{j\in[n]:i_j\in R_y\}|/n\le v+\varepsilon$ for every $y\in Y$.
\end{enumerate}
\end{fact}
\begin{proof}
By Fact~\ref{fact:minimax} there are distributions $\nu$ on $Y$ and $\mu$ on $I$ such that $\nu(E_i)\ge v$ for all $i\in I$ and $\mu(R_y)\le v$ for all $y\in Y$. Apply the $\varepsilon$-approximation theorem to $(\E,\nu)$ for (i) and to $(\E^*,\mu)$ for (ii) \cite[Theorem~2.1 and Lemma~2.2]{MY}.
\end{proof}

\begin{definition}\label{def:transversal-pq}
Let $\E\coloneqq(E_i)_{i\in I}$ be an indexed family of subsets of $Y$.
\begin{enumerate}[label=\textup{(\roman*)},leftmargin=*]
\item A set $T\subseteq Y$ is a \emph{transversal} of $\E$ if $T\cap E_i\ne\varnothing$ for every $i\in I$.
\item For integers $p\ge q\ge1$, the indexed family $\E$ has the \emph{$(p,q)$-property} if among any $p$ pairwise distinct sets $E_{i_1},\ldots,E_{i_p}$ some $q$ have a common point.
\end{enumerate}
\end{definition}

Matou\v{s}ek proved that if a family has dual shatter function $\pi^*(m)=o(m^q)$, then for each $p\ge q$ its finite subfamilies of nonempty sets with the $(p,q)$-property have uniformly bounded transversal size \cite[Theorem~4]{Mat}. He deduces his $(p,q)$-theorem from his fractional Helly theorem by following the proof of the $(p,q)$-theorem of Alon and Kleitman for families of convex sets \cite{AK} (their proof uses linear programming duality and weak $\varepsilon$-nets, whose roles here are supplied by Fact~\ref{fact:vc-minimax}). Lemma~\ref{lem:pq} records an explicit quantitative refinement.

\begin{fact}[{Matou\v{s}ek's fractional Helly bound \cite[Theorem~2]{Mat}}]\label{fact:frac-helly}
Let $\E\coloneqq(E_i)_{i\in I}$ be a finite indexed family of subsets of $Y$ with dual shatter function $\pi^*$. Let $1\le q\le t$ be integers with $t\ge2$, and suppose $0<\alpha\le1$ satisfies $4\pi^*(t)<\alpha\binom tq$. For $n\ge2t^2$, let $(i_1,\ldots,i_n)\in I^n$. If at least $\alpha\binom nq$ of the $q$-element subsets $J\subseteq[n]$ satisfy $\bigcap_{j\in J}E_{i_j}\ne\varnothing$, then some $y\in Y$ satisfies $|\{j\in[n]:y\in E_{i_j}\}|\ge n/(2t)$.
\end{fact}

\begin{proof}
Under the hypothesis $4\pi^*(t)<\alpha\binom tq$, Matou\v{s}ek's argument gives a point common to at least $\beta n$ entries for $\beta\coloneqq1/(2t)$, provided $n\beta\ge t$. This last condition is precisely $n\ge2t^2$. The argument applies to repeated indices because the traces on any $t$ positions of the tuple are determined by at most $t$ distinct indices.
\end{proof}

\begin{lemma}[A polynomial $(p,q)$-theorem]\label{lem:pq}
There is an absolute integer $\constof{lem:pq}\ge1$ with the following property. Let $\E\coloneqq(E_i)_{i\in I}$ be a nonempty finite indexed family of nonempty subsets of a finite set $Y$, with $\VC(\E)\le d$ and $\pi^*_{\E}(m)\le(m+1)^s$ for all $m\ge1$, where $d,s\ge1$. If $\E$ has the $(p,q)$-property for integers $p\ge q>s$, then it has a transversal of size $O_{d,s,q}(p^{2q/(q-s)})$, and of size at most $\constof{lem:pq}(d+1)p^4$ if $q\ge2s$.
\end{lemma}
\begin{proof}
Put $t\coloneqq\lceil(32p)^{q/(q-s)}\rceil$ and $\alpha\coloneqq q!/(2p)^q$. Since $t^{q-s}\ge(32p)^q$ and $t\ge2q$, we have
\[
 4\pi^*_{\E}(t)\le4(2t)^s<8^qt^s\le(t/(4p))^q\le\alpha\binom tq .
\]
Let $v$ be the value of the incidence matrix of $\E$. Take $(i_1,\ldots,i_n)\in I^n$ from Fact~\ref{fact:vc-minimax}(ii) for $\varepsilon\coloneqq1/(8t)$. Put $m\coloneqq nt^2$, $r\coloneqq2pt^2$, and $N\coloneqq nr=2pm\ge2t^2$. Define an indexed family $(F_j)_{j\in[N]}$ of subsets of $Y$ by $F_{(h-1)r+a}\coloneqq E_{i_h}$ for $h\in[n]$ and $a\in[r]$ (so each occurrence of a set in $(E_{i_1},\ldots,E_{i_n})$ is repeated $r$ times). For $H_y\coloneqq\{j\in[N]:y\in F_j\}$ we have
\[
 \frac{|H_y|}{N}
 =\frac{r\,|\{h\in[n]:y\in E_{i_h}\}|}{rn}\le v+\varepsilon
 \qquad\text{for all }y\in Y.
\]
Suppose $v+\varepsilon<1/(2t)$. We claim that there are $m$-element sets $B_1,\ldots,B_p\subseteq[N]$ such that the families $\{F_j:j\in B_b\}$, $b\in[p]$, are pairwise disjoint.
For each $A\in\{F_1,\ldots,F_N\}$ put $C_A\coloneqq\{j\in[N]:F_j=A\}$. Taking any $y\in A$ gives $C_A\subseteq H_y$, so $|C_A|<N/(2t)\le m$. To choose $B_b$, take unused sets $C_A$ until their union $U$ first contains at least $m$ indices. Since the last set added has size less than $m$, we have $m\le|U|<2m$. Choose any $m$ indices in $U$ for $B_b$ and discard all of $U$. Before step $b\le p$ at least $N-2m(b-1)\ge2m$ indices remain, so all $p$ steps are possible. Discarding $U$ removes every copy of each set represented in $B_b$, proving the claim.
Now for each $(j_1,\ldots,j_p)\in B_1\times\cdots\times B_p$, the $(p,q)$-property gives some $J\subseteq\{j_1,\ldots,j_p\}$ of size $q$ with $\bigcap_{j\in J}F_j\ne\varnothing$. Each such $J$ occurs in at most $m^{p-q}$ of the $m^p$ choices, since it specifies the choice of $j_b$ for $q$ indices $b$. Thus at least $m^q\ge\alpha\binom Nq$ of the $q$-subsets of $[N]$ satisfy this intersection condition. Fact~\ref{fact:frac-helly} gives $y\in Y$ with $|H_y|/N\ge 1/(2t)$, contradicting $|H_y|/N\le v+\varepsilon<1/(2t)$.

Consequently $v-\varepsilon\ge1/(2t)-2\varepsilon=1/(4t)$. Fact~\ref{fact:vc-minimax}(i) gives $(y_1,\ldots,y_k)\in Y^k$ with $k\le \constof{fact:vc-minimax}(d+1)(8t)^2$ and $|\{j\in[k]:y_j\in E_i\}|/k\ge v-\varepsilon>0$ for every $i\in I$. Hence $T\coloneqq\{y_j:j\in[k]\}$ meets every $E_i$. The bound on $k$ gives the claimed estimates, since $t\le(32p)^2$ when $q\ge2s$.
\end{proof}

\section{Construction of the scheme}\label{sec:descent}
\subsection{Proper reconstruction or disagreement location relative to a fixed labeling and subclass}\label{sec:correction}
\begin{definition}
A \emph{partial labeling} on $X$ is a map $h:A\to\bits$ with domain $\dom h\coloneqq A\subseteq X$.
\end{definition}

\begin{lemma}\label{lem:correct}
There is an integer $\constof{lem:correct}(d)=2^{O(d)}$ such that for every nonempty $\Z\subseteq\C$ and every partial labeling $h$ there are inclusion-preserving maps $R:2^X\to2^\Z$ and $P:2^X\to2^{\dom h}$ satisfying the following. If $B\subseteq\dom h$, $a\in\Z$ and $a|_B=h|_B$, then at least one of the following holds:
\begin{itemize}[nosep]
\item Some $c\in R(B)$ satisfies $c|_B=a|_B$.
\item Some $y\in P(B)$ satisfies $h(y)\ne a(y)$.
\end{itemize}
For every $T\subseteq X$, every $c\in R(T)$ belongs to $R(T_0)$ for some $T_0\subseteq T$ of size at most one, and every $y\in P(T)$ belongs to $P(T_0)$ for some $T_0\subseteq T$ of size at most four. Moreover,
\[
 |R(T)|\le \constof{lem:correct}(1+|T|),\qquad
 |P(T)|\le \constof{lem:correct}|T|^4\qquad\text{for all }T\subseteq X.
\]
\end{lemma}
\begin{proof}
If some $c\in\Z$ extends $h$, put $R(T)\coloneqq\{c\}$ and $P(T)\coloneqq\varnothing$ for every $T$. Otherwise all the sets $E_c\coloneqq\{y\in\dom h:c(y)\ne h(y)\}$, $c\in\Z$, are nonempty. Set $D\coloneqq2^{d+1}-1$ and $K\coloneqq2D-1$. Consider tuples $(c_1,\ldots,c_m)\in\Z^m$ such that $|\{i\in[m]:y\in E_{c_i}\}|\le K$ for every $y\in\dom h$. Since each $E_{c_i}$ is nonempty, $m\le\sum_{i=1}^m |E_{c_i}|=\sum_{y\in\dom h}|\{i\in[m]:y\in E_{c_i}\}|\le K|\dom h|$, so the possible lengths are bounded. Taking $K$ copies of any concept in $\Z$ gives a tuple of length $K$. We may therefore fix a tuple $(c_1,\ldots,c_L)$ of maximum length $L\ge K$.

By maximality of $L$, among any $L+1$ distinct sets in $\E\coloneqq(E_c)_{c\in\Z}$ some $K+1$ have a common point, so $\E$ has the $(L+1,K+1)$-property. By Lemma~\ref{lem:vc-bounds}(i)--(ii), $\VC(\E)\le d$ and $\pi^*_\E(t)\le(t+1)^D$. Since $K+1=2D$, Lemma~\ref{lem:pq} gives a transversal $V\subseteq\dom h$ of $\E$ with
\[
 |V|\le \constof{lem:pq}(d+1)(L+1)^4\le A L^4,\qquad
 A\coloneqq16(d+1)\constof{lem:pq}.
\]
Thus we may choose $\gamma:V\to[L]^4$ with fibers of size at most $A$.

For $T\subseteq X$, let $Q(T)\coloneqq\{i\in[L]:T\cap E_{c_i}\ne\varnothing\}$. Let
\[
 R(T)\coloneqq\{c_1\}\cup\{c_{i+1}:i\in Q(T),\ i<L\},\qquad
 P(T)\coloneqq\{y\in V:\gamma(y)\in Q(T)^4\}.
\]
Both sets increase with $T$ since $Q(T)$ does. We have $c_1\in R(\varnothing)$. For any $c\in R(T)\setminus\{c_1\}$, choose $i\in Q(T)$ with $i<L$ and $c=c_{i+1}$. Take $x\in T\cap E_{c_i}$. Then $i\in Q(\{x\})$, so $c\in R(\{x\})$.

For $y\in P(T)$, put $(i_1,i_2,i_3,i_4)\coloneqq\gamma(y)\in Q(T)^4$. For each $j\in[4]$, choose $x_j\in T\cap E_{c_{i_j}}$ and put $T_0\coloneqq\{x_1,x_2,x_3,x_4\}$. Then $T_0\subseteq T$, $|T_0|\le4$, and $i_j\in Q(T_0)$ for every $j$. Thus $\gamma(y)\in Q(T_0)^4$ and $y\in P(T_0)$.

For each $i\in Q(T)$ there exists $y\in T\cap E_{c_i}$, and each $y\in T\cap\dom h$ belongs to at most $K$ of the $E_{c_i}$, counting indices. Thus $|Q(T)|\le\sum_{y\in T\cap\dom h}|\{i\in[L]:y\in E_{c_i}\}|\le K|T|$. By definition $|R(T)|\le1+|Q(T)|$ and, since each fiber of $\gamma$ has size at most $A$, $|P(T)|\le A|Q(T)|^4$. The cardinality bounds on $R(T)$ and $P(T)$ hold with $\constof{lem:correct}(d)\coloneqq\max\{K+1,AK^4\}$.

If some $c_i$ agrees with $a$ on $B$, take the least such $i$. Either $i=1$, or $c_{i-1}$ disagrees with $a|_B=h|_B$ at a point of $B$, so $i-1\in Q(B)$. In both cases $c_i\in R(B)$. If no $c_i$ agrees with $a$ on $B$, then $Q(B)=[L]$ and $P(B)=V$. Since $V$ meets $E_a$, some $y\in P(B)$ satisfies $h(y)\ne a(y)$.
\end{proof}

\begin{definition}\label{def:fixed-maps}
For each nonempty subclass $\Z\subseteq\C$ and partial labeling $h$, fix one pair of maps constructed in the proof of Lemma~\ref{lem:correct}. Denote their values at $T$ by $R_{\Z,h}(T)$ and $P_{\Z,h}(T)$.
\end{definition}

\subsection{Partial labelings determined by sampled points}
The dimension descent will choose a tentative labeling from a tuple of shattered subsets of the sample. Each subset is accompanied in the proof by concepts witnessing that it is shattered. The next lemma bounds how many such tuples of concepts suffice.

\begin{lemma}\label{lem:maj}
For each integer $d\ge1$ there is a positive integer $\constof{lem:maj}(d)=O(d2^d)$ with the following property. Let $\C\subseteq\bits^X$ have $\VC(\C)\le d$. Fix $a\in\C$, $Y\subseteq X$, and a nonempty finite set $I$. Put $m\coloneqq2^d$. For each $i\in I$, let $m_i\in[m]$ and $w_i\coloneqq(c_{i,j})_{j\in[m_i]}\in\C^{m_i}$, and define
\[
 E_i\coloneqq\{y\in Y:c_{i,j}(y)\ne a(y)\text{ for some }j\in[m_i]\}.
\]
Put $\E\coloneqq(E_i)_{i\in I}$. Then either $\E$ has a transversal of size at most $d\,\constof{lem:maj}$, or there is a sequence $(i_1,\ldots,i_u)\in I^u$ of length $1\le u\le\constof{lem:maj}$ such that $|\{j\in[u]:y\in E_{i_j}\}|<u/2$ for every $y\in Y$.
\end{lemma}
\begin{proof}
If $Y=\varnothing$, take $u\coloneqq1$ and any $i_1\in I$. Otherwise let $\mathcal D\coloneqq(D_c)_{c\in\C}$, where $D_c\coloneqq\{y\in Y:c(y)\ne a(y)\}$. Lemma~\ref{lem:vc-bounds}(i)--(ii) gives $\VC(\mathcal D)\le d$ and $\VC(\mathcal D^*)\le 2^{d+1}-1$. Each $E_i=\bigcup_{j=1}^{m_i}D_{c_{i,j}}$ belongs to the family of all $m$-fold unions from $\mathcal D$, since its last set can be repeated $m-m_i$ times. Lemma~\ref{lem:vc-bounds}(iii) therefore gives $\VC(\E)=O(d^2 2^d)$ and $\VC(\E^*)=O(d2^d)$.

Apply Fact~\ref{fact:vc-minimax} to $\E$ with $\varepsilon\coloneqq1/8$. If the incidence matrix of $\E$ has value $v\ge1/4$, each $E_i$ contains at least $1/8$ of the points in \textup{(i)}, so these points form a transversal. Otherwise \textup{(ii)} gives $|\{j\in[u]:y\in E_{i_j}\}|/u<3/8$ for every $y\in Y$. Since $\varepsilon^{-2}=64$, part~(i) uses at most $64\constof{fact:vc-minimax}(\VC(\E)+1)$ points and part~(ii) at most $64\constof{fact:vc-minimax}(\VC(\E^*)+1)$ indices. By the preceding dimension estimates these lengths are at most $cd^2 2^d$ and $cd2^d$, respectively, for an absolute $c\ge2$. Take $\constof{lem:maj}(d)\coloneqq\lceil cd2^d\rceil$.
\end{proof}

Take $\constof{lem:maj}(d)$ as in Lemma~\ref{lem:maj}. Write $\binom Xr$ for the set of $r$-element subsets of $X$.

\begin{definition}\label{def:supported-labeling}
Fix $\varnothing\ne\Z\subseteq\C$, $r\in[d]$, $u\in[\constof{lem:maj}]$ and $\mathbf S\coloneqq(S_1,\ldots,S_u)\in(\binom Xr)^u$.
\begin{enumerate}[label=\textup{(\roman*)},leftmargin=*,nosep]
\item For $y\in X$ and $e\in\bits$ let $\Z_{y,e}\coloneqq\{c\in\Z:c(y)=e\}$ and $I_e(y)\coloneqq\{i\in[u]:\Z_{y,e}\text{ shatters }S_i\}$. The bit $e$ is \emph{supported at $y$} if $y\in U_e$, where $U_e\coloneqq\{y\in X:|I_e(y)|>u/2\}$.
\item The partial labeling associated to $(\Z,\mathbf S)$ is
\[
 h_{\Z,\mathbf S}:X\setminus(U_0\cap U_1)\longrightarrow\bits,\qquad
 h_{\Z,\mathbf S}(y)\coloneqq\mathbf 1_{U_1}(y).
\]
In particular $h_{\Z,\mathbf S}(y)=0$ when $y\notin U_0\cup U_1$.
\item For $S\subseteq X$, a \emph{witness} for shattering $S$ by $\Z$ is a tuple $(c_\eta)_{\eta\in\bits^S}$ with $c_\eta\in\Z$ and $c_\eta|_S=\eta$ for every $\eta\in\bits^S$.
\end{enumerate}
\end{definition}

\begin{lemma}\label{lem:support}
Let $\Z,r,u,\mathbf S$ be as in Definition~\ref{def:supported-labeling}. Suppose $B\subseteq X$, $\VC(\Z|_B)=r$ and $S_i\subseteq B$ for all $i\in[u]$. Then $B\subseteq\dom h$ for $h\coloneqq h_{\Z,\mathbf S}$. Suppose further that the tuple $(c_{i,\eta})_{\eta\in\bits^{S_i}}$ witnesses shattering of $S_i$ by $\Z$ for each $i$. For $y\in\dom h$ and $e\in\bits$, we have $h(y)=e$ whenever
\begin{equation}\label{eq:witness-majority}
 \bigl|\{i\in[u]:c_{i,\eta}(y)=e\ \text{for every }\eta\in\bits^{S_i}\}\bigr|>u/2.
\end{equation}
\end{lemma}
\begin{proof}
Suppose some $b\in B$ belongs to $U_0\cap U_1$, i.e., both bits are supported at $b$. Since $|I_0(b)|+|I_1(b)|>u$, some $i\in[u]$ belongs to both $I_0(b)$ and $I_1(b)$. Thus both $\Z_{b,0}$ and $\Z_{b,1}$ shatter $S_i$. In particular $b\notin S_i$, since each subclass is constant at $b$. Fix any labeling $\eta:S_i\cup\{b\}\to\bits$ and put $e\coloneqq\eta(b)$. Since $\Z_{b,e}$ shatters $S_i$, some $c\in\Z_{b,e}$ satisfies $c|_{S_i}=\eta|_{S_i}$. By definition $c(b)=e=\eta(b)$, so $c|_{S_i\cup\{b\}}=\eta$. Hence $\Z$ shatters the $(r+1)$-element subset $S_i\cup\{b\}$ of $B$, contradicting $\VC(\Z|_B)=r$. Thus $B\subseteq\dom h$.

For the second assertion, each index $i\in[u]$ counted in \eqref{eq:witness-majority} belongs to $I_e(y)$, since for every $\eta\in\bits^{S_i}$ the concept $c_{i,\eta}$ belongs to $\Z_{y,e}$ and restricts to $\eta$ on $S_i$. Hence $|I_e(y)|>u/2$, so $e$ is supported at $y$. Since $y\in\dom h$, it follows that $h(y)=e$.
\end{proof}

\begin{definition}\label{def:descriptions}
Fix a linear order on $\C$. For $\varnothing\ne\Z\subseteq\C$ and $T\subseteq X$, let $\mathcal H_\Z(T)$ consist of the partial labelings $h_{\Z,\mathbf S}$ from Definition~\ref{def:supported-labeling} for which every $S_i$ is contained in $T$. Here all $r\in[d]$ and $1\le u\le\constof{lem:maj}$ are allowed. Define the \emph{available points} and \emph{available concepts} in $\Z$ from $T$ by
\[
 P_\Z(T)\coloneqq T\cup\!\!\bigcup_{h\in\mathcal H_\Z(T)} P_{\Z,h}(T),
 \qquad
 R_\Z(T)\coloneqq\{\min\Z\}\cup\!\!\bigcup_{h\in\mathcal H_\Z(T)} R_{\Z,h}(T).
\]
\end{definition}

\begin{lemma}\label{lem:availability}
There is a positive integer $\constof{lem:availability}(d)=O(d^2 2^d)$ with the following properties. For every nonempty $\Z\subseteq\C$, both $P_\Z(T)$ and $R_\Z(T)$ increase with $T\subseteq X$. For each $y\in P_\Z(T)$ there is a subset $T_0\subseteq T$ of size at most $\constof{lem:availability}$ with $y\in P_\Z(T_0)$. Similarly, for each $c\in R_\Z(T)$ there is a subset $T_0\subseteq T$ of size at most $\constof{lem:availability}$ with $c\in R_\Z(T_0)$. For $t\coloneqq|T|$,
\[
 \max\{|P_\Z(T)|,|R_\Z(T)|\}\le M(t)\coloneqq(t+2)^{\constof{lem:availability}}
   =(t+2)^{O(d^2 2^d)}.
\]
\end{lemma}
\begin{proof}
If $T\subseteq T'$, then $\mathcal H_\Z(T)\subseteq\mathcal H_\Z(T')$, since each tuple of subsets of $T$ is also a tuple of subsets of $T'$. Monotonicity now follows from Definition~\ref{def:descriptions} and the inclusion-preserving maps of Lemma~\ref{lem:correct}.

Put $m\coloneqq d\,\constof{lem:maj}$. For $y\in T$ take $T_0\coloneqq\{y\}$, and for $c=\min\Z$ take $T_0\coloneqq\varnothing$. For the remaining cases, fix $h\in\mathcal H_\Z(T)$, choose a tuple $\mathbf S$ of subsets of $T$ with $h_{\Z,\mathbf S}=h$, and put $U\coloneqq\bigcup_{i=1}^u S_i$. Then $|U|\le ru\le m$ and $h\in\mathcal H_\Z(U)$. For $y\in P_{\Z,h}(T)$, Lemma~\ref{lem:correct} gives $V\subseteq T$ of size at most four with $y\in P_{\Z,h}(V)$. Taking $T_0\coloneqq U\cup V$ gives $h\in\mathcal H_\Z(T_0)$ and $y\in P_{\Z,h}(V)\subseteq P_{\Z,h}(T_0)\subseteq P_\Z(T_0)$. For $c\in R_{\Z,h}(T)$ the same argument uses $V$ of size at most one. Thus $|T_0|\le m+4$ in both cases.

There are at most $m$ choices of $(r,u)$ and at most $\binom tr^u\le(t+1)^{ru}$ tuples $\mathbf S$ for each choice, so $|\mathcal H_\Z(T)|\le m(t+1)^m$. By Lemma~\ref{lem:correct}, each $h$ contributes at most $\constof{lem:correct}(t+1)^4$ elements to either union. Including $T$ or $\{\min\Z\}$ gives a bound of $(1+m\constof{lem:correct})(t+1)^{m+4}$ for both cardinalities. Since $t+2\ge2$, this is at most $(t+2)^{\constof{lem:availability}}$ for $\constof{lem:availability}(d)\coloneqq m+4+\lceil\log_2(1+m\constof{lem:correct})\rceil$. This constant also bounds $|T_0|$. It is $O(d^2 2^d)$ since $m=O(d^2 2^d)$ and $\log_2\constof{lem:correct}=O(d)$.
\end{proof}

\subsection{Proper reconstruction or dimension descent}\label{sec:dimension-descent}
\begin{lemma}\label{lem:dimension-descent}
Let $a\in\Z\subseteq\C$, $B\subseteq X$, and $r\coloneqq\VC(\Z|_B)>0$. At least one of the following holds.
\begin{enumerate}[label=\textup{(\roman*)},leftmargin=*,nosep]
\item Some $c\in R_\Z(B)$ satisfies $c|_B=a|_B$.
\item There is a set $F\subseteq P_\Z(B)$ of at most $\constof{lem:availability}$ points such that $\Z_F\coloneqq\{c\in\Z:c|_F=a|_F\}$ has $\VC(\Z_F|_B)<r$.
\end{enumerate}
\end{lemma}
\begin{proof}
Put $Y\coloneqq P_\Z(B)$, so $B\subseteq Y$. Let $I$ consist of all pairs $(S,w)$ where $S\subseteq B$, $|S|=r$, and $w=(c_\eta)_{\eta\in\bits^S}$ with $c_\eta\in\Z$ and $c_\eta|_S=\eta$ for every $\eta\in\bits^S$. For $(S,w)\in I$ let $E_{(S,w)}\coloneqq\{y\in Y:c_\eta(y)\ne a(y)\text{ for some }\eta\in\bits^S\}$. For $F\subseteq Y$ put $\Z_F\coloneqq\{c\in\Z:c|_F=a|_F\}$. Then
\begin{equation}\label{eq:transversal-dimension}
 F\text{ meets every }E_{(S,w)}
 \quad\Longleftrightarrow\quad \VC(\Z_F|_B)<r .
\end{equation}
If $F\cap E_{(S,w)}=\varnothing$, then $c_\eta|_F=a|_F$ for every $\eta\in\bits^S$. Thus each $c_\eta$ belongs to $\Z_F$ and agrees with $\eta$ on $S$, so $\Z_F$ shatters $S$ (and $|S|=r$). Conversely, if $\Z_F$ shatters an $r$-element $S\subseteq B$, choose $c_\eta\in\Z_F$ with $c_\eta|_S=\eta$ for each $\eta\in\bits^S$. Their tuple $w$ satisfies $F\cap E_{(S,w)}=\varnothing$ because every $c_\eta$ agrees with $a$ on $F$. Since $\VC(\Z_F|_B)\le r$, this proves \eqref{eq:transversal-dimension}.

Since $\VC(\Z|_B)=r$, there is an $r$-element $S\subseteq B$ with $\Z|_S=\bits^S$. Choosing one extension in $\Z$ of each labeling of $S$ gives $(S,w)\in I$. For every $(S,w)\in I$, the tuple $w$ has $|\bits^S|=2^r\le 2^d$ entries. If Lemma~\ref{lem:maj}, applied on $Y$ with target $a$, gives a transversal $F$ of the indexed family $(E_{(S,w)})_{(S,w)\in I}$ with $|F|\le d\,\constof{lem:maj}\le\constof{lem:availability}$, then \eqref{eq:transversal-dimension} gives~\textup{(ii)}.

Otherwise let $((S_j,w_j))_{j=1}^u\in I^u$, with $u\in[\constof{lem:maj}]$, be the sequence of indices selected by Lemma~\ref{lem:maj}. For $j\in[u]$ and $\eta\in\bits^{S_j}$, denote the $\eta$-entry of $w_j$ by $c_{j,\eta}$. For every $y\in Y$,
\[
 \bigl|\{j\in[u]:c_{j,\eta}(y)=a(y)\text{ for every }\eta\in\bits^{S_j}\}\bigr|
 =u-\bigl|\{j\in[u]:y\in E_{(S_j,w_j)}\}\bigr|>u/2.
\]
Thus \eqref{eq:witness-majority} holds with $a(y)$ in place of $e$. Put $h\coloneqq h_{\Z,(S_1,\ldots,S_u)}$. It follows by Lemma~\ref{lem:support} that $B\subseteq\dom h$ and $h(y)=a(y)$ for all $y\in Y\cap\dom h$. Since $B\subseteq Y$ and $h\in\mathcal H_\Z(B)$, both $B$ and $P_{\Z,h}(B)$ are contained in $Y\cap\dom h$. Hence $h|_B=a|_B$ and no $y\in P_{\Z,h}(B)$ satisfies $h(y)\ne a(y)$. Lemma~\ref{lem:correct} therefore gives a concept $c\in R_{\Z,h}(B)\subseteq R_\Z(B)$ with $c|_B=a|_B$, proving~\textup{(i)}.
\end{proof}

\subsection{The compression and reconstruction maps}\label{sec:maps}
\begin{definition}\label{def:candidate-tree}
For $T\subseteq X$, define a rooted tree $\mathcal T(T)$ recursively. Its root is the one-term sequence $(\C)$. A node $v=(\Z_0,\ldots,\Z_j)$ has depth $j$ and associated subclass $\Z(v)\coloneqq\Z_j\subseteq\C$. If $j<d$, then for each $F\subseteq P_{\Z_j}(T)$ with $|F|\le\constof{lem:availability}$ and each labeling $\eta\in\Z_j|_F$, put
\[
 \Z'\coloneqq\{c\in\Z_j:c|_F=\eta\}.
\]
The children of $v$ are precisely the sequences $(\Z_0,\ldots,\Z_j,\Z')$ obtained in this way. Each $\Z'$ is nonempty since $\eta\in\Z_j|_F$. Deleting the last subclass gives the parent of each nonroot node. Nodes of depth $d$ have no children.

Let $L(T)$ be the list of concepts in $\bigcup_{v\in\mathcal T(T)}R_{\Z(v)}(T)$, each appearing once in the fixed linear order on $\C$.
\end{definition}

\begin{lemma}\label{lem:retained-set}
Put $k\coloneqq(d\constof{lem:availability}+1)\constof{lem:availability}$. For every realizable sample $(B,\ell)$ there exist $T\subseteq B$ with $|T|\le k$ and a concept $c\in L(T)$ with $c|_B=\ell$.
\end{lemma}
\begin{proof}
By realizability, choose $a\in\C$ with $a|_B=\ell$. Choose a node $(\Z_0,\ldots,\Z_j)$ of maximum depth in $\mathcal T(B)$ subject to $a\in\Z_j$ and $\VC(\Z_{i+1}|_B)<\VC(\Z_i|_B)$ for every $0\le i<j$. The root $(\C)$ satisfies these conditions, and all depths are at most $d$.

We claim that some $c\in R_{\Z_j}(B)$ extends $\ell$. Put $r\coloneqq\VC(\Z_j|_B)$. If $r=0$, every member of $\Z_j$ agrees with $a$ on $B$, so take $c\coloneqq\min\Z_j$. If $r>0$, the $j$ strict decreases give $j+r\le\VC(\C|_B)\le d$, hence $j<d$. If no such $c$ existed, Lemma~\ref{lem:dimension-descent} would give $F\subseteq P_{\Z_j}(B)$ of size at most $\constof{lem:availability}$ such that $\Z'\coloneqq\{c'\in\Z_j:c'|_F=a|_F\}$ has $\VC(\Z'|_B)<r$. Since $a\in\Z'$, appending $\Z'$ to the chosen node would give a child satisfying the same conditions, contrary to maximality of $j$.

For each $0\le i<j$, Definition~\ref{def:candidate-tree} gives $F_i\subseteq P_{\Z_i}(B)$ of size at most $\constof{lem:availability}$ with $\Z_{i+1}=\{c'\in\Z_i:c'|_{F_i}=a|_{F_i}\}$, since $a\in\Z_j\subseteq\Z_{i+1}$. For each $y\in F_i$, Lemma~\ref{lem:availability} gives $U_{i,y}\subseteq B$ of size at most $\constof{lem:availability}$ with $y\in P_{\Z_i}(U_{i,y})$. It also gives $U\subseteq B$ of size at most $\constof{lem:availability}$ with $c\in R_{\Z_j}(U)$. Let $T$ be the union of $U$ and all the $U_{i,y}$. There are at most $j\constof{lem:availability}\le d\constof{lem:availability}$ pairs $(i,y)$, so $|T|\le k$.

By monotonicity, $F_i\subseteq P_{\Z_i}(T)$ for every $0\le i<j$. Since $\Z_{i+1}=\{c'\in\Z_i:c'|_{F_i}=a|_{F_i}\}$ for each such $i$, Definition~\ref{def:candidate-tree} gives $(\Z_0,\ldots,\Z_j)\in\mathcal T(T)$. Thus $c\in R_{\Z_j}(U)\subseteq R_{\Z_j}(T)\subseteq L(T)$.
\end{proof}

\begin{lemma}\label{lem:candidate-count}
For every $T\subseteq X$ with $|T|\le k$, $|L(T)|\le N\coloneqq(d+1)M(k)(2M(k)+1)^{d\constof{lem:availability}}$. In particular, $k+\log_2 N=2^{O(d)}$.
\end{lemma}
\begin{proof}
Fix a node $v\in\mathcal T(T)$ and put $\Z\coloneqq\Z(v)$. By Lemma~\ref{lem:availability}, $|P_\Z(T)|,|R_\Z(T)|\le M(k)$. By Definition~\ref{def:candidate-tree}, each child of $v$ is specified by a set $F\subseteq P_\Z(T)$ of size at most $\constof{lem:availability}$ and a labeling $\eta\in\Z|_F$. Write the pairs $(y,\eta(y))$ for $y\in F$ in any order and pad with a blank symbol to length $\constof{lem:availability}$. There are at most $2M(k)+1$ choices per position, so each node has at most $(2M(k)+1)^{\constof{lem:availability}}$ children. Each of the $d+1$ levels therefore has at most $(2M(k)+1)^{d\constof{lem:availability}}$ nodes. Since each node contributes at most $M(k)$ concepts to $L(T)$, we obtain $|L(T)|\le N$.

By Lemma~\ref{lem:availability}, $\constof{lem:availability}=O(d^2 2^d)$, so $k=O(d^5 4^d)$. Substituting $\log_2 M(k)=\constof{lem:availability}\log_2(k+2)$ in the formula for $N$ gives $\log_2 N=O(d^6 4^d)$.
\end{proof}

\begin{definition}[Compression map]\label{def:compression}
Take $k$ and $N$ from Lemmas~\ref{lem:retained-set} and~\ref{lem:candidate-count}, respectively. Put $I_d\coloneqq[N]$ and fix a linear order on the subsets of $X$. For $(B,\ell)\in L_\C(\infty)$, let $T$ be the first subset of $B$ of size at most $k$ for which $L(T)$ contains a concept extending $\ell$. Let $i$ be the first position of such a concept in $L(T)$. Define
\[
 \kappa:L_\C(\infty)\longrightarrow L_\C(k)\times I_d,\qquad
 \kappa(B,\ell)\coloneqq((T,\ell|_T),i).
\]
\end{definition}
\begin{definition}[Reconstruction map]\label{def:reconstruction}
For $((T,\lambda),i)\in L_\C(k)\times I_d$ let $\rho((T,\lambda),i)$ be the $i$th concept of $L(T)$ if $i\le|L(T)|$, and $\min\C$ otherwise.
\end{definition}

By its definition, $\rho$ takes values in $\C$, and the definition of $\kappa$ gives $\rho(\kappa(B,\ell))|_B=\ell$. This proves Theorem~\ref{thm:main} for finite $X$ and $d\ge1$. (For $d=0$ the class has a single member. The encoder always retains the empty sample with a single fixed tag, and the decoder always returns that member.)

\begin{remark}
The reconstruction map in Definition~\ref{def:reconstruction} uses only $T$ and $i$. On an encoded sample, $i$ selects a concept $c\in L(T)$ whose restriction to $T$ is the retained labeling $\lambda$. More generally, any binary scheme retaining at most $k$ examples can include $\lambda$ in its tag using $k$ additional bits.
\end{remark}

\subsection{Arbitrary domains}\label{sec:arbitrary}
Let $X$ be infinite and let $\varnothing\ne\C\subseteq\bits^X$ have $\VC(\C)\le d$, where $d\ge1$. Let $M$ be the structure on $X$ with a unary predicate $U_A$ naming each $A\subseteq X$. For each integer $n\ge1$ and family $\mathcal F\subseteq 2^X$ whose finite subfamilies have transversals of size at most $n$, introduce a separate tuple of constants $z_1,\ldots,z_n$ and impose $\bigvee_{j=1}^n U_E(z_j)$ for every $E\in\mathcal F$. Each finite collection of these conditions is satisfiable in $M$ by the transversal bounds. Compactness gives an elementary extension $\widehat M$ satisfying them all. Write $\widehat X$ for its universe and $\widehat A\coloneqq U_A(\widehat M)$ for $A\subseteq X$. For $h:A\to\bits$, let $\widehat h:\widehat A\to\bits$ have fibers $\widehat{h^{-1}(e)}$ for $e\in\bits$. These fibers partition $\widehat A$ by elementarity.

\paragraph{Extension of Lemma~\ref{lem:correct}.}
Fix $\varnothing\ne\Z\subseteq\C$ and $h:A\to\bits$, where $A\subseteq X$. We construct the maps of Lemma~\ref{lem:correct} for finite $T\subseteq X$, still with $R(T)\subseteq\Z$ but now with $P(T)\subseteq\widehat A$. For finite $B\subseteq A$ and $a\in\Z$ with $a|_B=h|_B$, we need either $c|_B=a|_B$ for some $c\in R(B)$ or $\widehat h(y)\ne\widehat a(y)$ for some $y\in P(B)$ (so $y$ may lie in $\widehat X\setminus X$).

For $c\in\Z$ put $E_c\coloneqq\{x\in A:c(x)\ne h(x)\}$. If $E_c=\varnothing$ for some $c$, take $R(T)\coloneqq\{c\}$ and $P(T)\coloneqq\varnothing$. Otherwise let $\E\coloneqq\{E_c:c\in\Z\}$, $D\coloneqq2^{d+1}-1$ and $K\coloneqq2D-1$. By Lemma~\ref{lem:vc-bounds}(i)--(ii) on finite restrictions, $\VC(\E)\le d$ and $\pi^*_\E(t)\le(t+1)^D$. Consider tuples $(c_1,\ldots,c_m)\in\Z^m$, $m\in\{1,2,\ldots\}$, satisfying $|\{i\in[m]:x\in E_{c_i}\}|\le K$ for every $x\in A$.

If $m$ is bounded, choose a tuple $(c_1,\ldots,c_L)$ of maximum length $L$. As in Lemma~\ref{lem:correct}, $L\ge K$ and $\E$ has the $(L+1,K+1)$-property. For a finite subfamily $\E_0\subseteq\E$, restrict $A$ to one point $x$ for each distinct set $\{E\in\E_0:x\in E\}$. Lemma~\ref{lem:pq} then gives a transversal of $\E_0$ of size at most $n\coloneqq\constof{lem:pq}(d+1)(L+1)^4$. By our choice of $\widehat M$, the family $\{\widehat E_c:c\in\Z\}$ has a transversal $W\subseteq\widehat X$ of size at most $n$. By elementarity, $\widehat E_c=\{y\in\widehat A:\widehat c(y)\ne\widehat h(y)\}$ for each $c\in\Z$, so $V\coloneqq W\cap\widehat A$ is a transversal in $\widehat A$ with $|V|\le n$. For each finite $T\subseteq X$, put $Q(T)\coloneqq\{i\in[L]:T\cap E_{c_i}\ne\varnothing\}$ and $R(T)\coloneqq\{c_1\}\cup\{c_{i+1}:i\in Q(T),\ i<L\}\subseteq\Z$. As in the proof of Lemma~\ref{lem:correct}, the bound $|V|\le n$ lets us choose $\gamma:V\to[L]^4$ with the same bound on its fiber sizes. Put $P(T)\coloneqq\gamma^{-1}(Q(T)^4)\subseteq V\subseteq\widehat A$. The bounds of Lemma~\ref{lem:correct} follow as before. If no $c_i$ agrees with $a$ on $B$, then $Q(B)=[L]$ and hence $P(B)=V$, and every $y\in V\cap\widehat E_a$ satisfies $\widehat h(y)\ne\widehat a(y)$.

If $m$ is unbounded, put $q\coloneqq K+1=2D$. Form the hypergraph with vertex set $\E$ whose edges are the $q$-element subfamilies with nonempty intersection. A tuple of length $m$ contains at least $m/K$ distinct sets, since each $E_c$ is nonempty. No $q$ of those sets have a common point, so the hypergraph has arbitrarily large finite independent sets. Since $\pi^*_\E(t)=o(t^q)$, by Matou\v{s}ek's fractional Helly theorem (Fact~\ref{fact:frac-helly}), for each $0<\alpha\le1$ there is a $\beta>0$ with the following property. If at least $\alpha\binom{m}{q}$ of the $q$-element subfamilies of $m\ge q$ distinct members of $\E$ have nonempty intersection, at least $\beta m$ of the $m$ sets contain a common point. Those members form a clique. Thus the infinitary $(p,q)$-theorem of Jung and P\'alv\"olgyi \cite[Theorem~4]{JP} gives an infinite independent set, from which we choose distinct $E_{c_1},E_{c_2},\ldots$. Set $P(T)\coloneqq\varnothing$ and
\[
 R(T)\coloneqq\{c_1\}\cup\{c_{i+1}: i\ge1,\ T\cap E_{c_i}\ne\varnothing\}.
\]
Each $x\in T$ belongs to at most $K$ of the $E_{c_i}$, so $|R(T)|\le1+K|T|$ and every member other than $c_1$ belongs to $R(\{x\})$ for some $x\in T$. A finite $B\subseteq A$ meets at most $K|B|$ of the $E_{c_i}$, so there is a least $i$ with $B\cap E_{c_i}=\varnothing$. Then $c_i|_B=h|_B$, and $c_i\in R(B)$ since either $i=1$ or $B\cap E_{c_{i-1}}\ne\varnothing$.

\paragraph{Dimension descent and reconstruction.}
For each $\Z,h$, fix a pair of maps $R_{\Z,h}$ and $P_{\Z,h}$ from the preceding construction. Use these maps in Definition~\ref{def:descriptions}, taking $\min\Z$ in a fixed well-order on $\C$. The proof of Lemma~\ref{lem:availability} uses only finite $T\subseteq X$ and the bounds of Lemma~\ref{lem:correct}, so it gives the same bounds for $P_\Z(T)\subseteq\widehat X$ and $R_\Z(T)\subseteq\Z$.

Fix $a\in\Z$, a finite $B\subseteq X$, and $r\coloneqq\VC(\Z|_B)>0$. Put $Y\coloneqq P_\Z(B)$, a finite subset of $\widehat X$ containing $B$. In the proof of Lemma~\ref{lem:dimension-descent}, choose one witness pair $(S,w)$ for each distinct $E_{(S,w)}\subseteq Y$, evaluating the concepts at points of $Y$ by their extensions. The finite class $\{\widehat c|_Y:c\in\Z\}$ has VC dimension at most $d$. Indeed, for any fixed $2^{d+1}$ concepts in $\Z$, the existence of $d+1$ points on which they realize all labelings is first-order in their fiber predicates, so it would transfer from $\widehat M$ to $M$. We may therefore apply Lemma~\ref{lem:maj} to these pairs in this class with target $\widehat a|_Y$. If Lemma~\ref{lem:maj} gives a transversal $F\subseteq Y\subseteq\widehat X$, it meets every $E_{(S,w)}$, so \eqref{eq:transversal-dimension} gives $\VC(\Z_F|_B)<r$ for $\Z_F\coloneqq\{c\in\Z:\widehat c|_F=\widehat a|_F\}$.

In the other case let $h\coloneqq h_{\Z,(S_1,\ldots,S_u)}$ on $X$ for the selected pairs. Since each $S_i\subseteq B$ has size $r$, the proof of Lemma~\ref{lem:support} gives $B\subseteq\dom h$. For fixed witnesses and $e\in\bits$, the implication in Lemma~\ref{lem:support} from \eqref{eq:witness-majority} to $h(y)=e$ on $\dom h$ is first-order in their fiber predicates and those of $h$. It holds in $\widehat M$, where the selected witnesses satisfy \eqref{eq:witness-majority} with $e=\widehat a(y)$ at each $y\in Y$. Hence $\widehat h=\widehat a$ on $Y\cap\widehat{\dom h}$. This intersection contains $B$ and, since $h\in\mathcal H_\Z(B)$, also $P_{\Z,h}(B)$. The extended Lemma~\ref{lem:correct} gives $c\in R_{\Z,h}(B)\subseteq R_\Z(B)$ with $c|_B=a|_B$.

In the tree of Definition~\ref{def:candidate-tree}, let $\Z(v)$ be the subclass associated to a node $v$. The child $v'$ determined by $F\subseteq P_{\Z(v)}(T)$ and $\eta\in\{\widehat c|_F:c\in\Z(v)\}$ has associated subclass $\Z(v')=\{c\in\Z(v):\widehat c|_F=\eta\}$. Thus every entry in $L(T)$ belongs to $\C\subseteq\bits^X$.

For each $y\in P_{\Z(v)}(B)$ used along a path, Lemma~\ref{lem:availability} gives $U\subseteq B$ of size at most $\constof{lem:availability}$ with $y\in P_{\Z(v)}(U)$. Hence for every realizable sample $(B,\ell)$, Lemmas~\ref{lem:retained-set} and~\ref{lem:candidate-count} give $T\subseteq B\subseteq X$ with $|T|\le k$ such that some $c$ extending $\ell$ occurs at a position $i\le N$ in $L(T)$. The maps of Definitions~\ref{def:compression} and~\ref{def:reconstruction} therefore encode $(B,\ell)$ as $((T,\ell|_T),i)$ and reconstruct the selected concept $c\in\C$, proving Theorem~\ref{thm:main}.

\begin{problem}
We have not tried to optimize the bound in Theorem~\ref{thm:main}. Does every class of VC-dimension $d$ admit a proper labeled compression scheme of total size polynomial in $d$, or even of size $O(d)$?
\end{problem}

\section*{Acknowledgements and AI Methodology}

When I learned of the question from Ben-David et al.~\cite[Section~3]{BDHMSY} in spring 2023, the distinction between proper and improper compression reminded me of that between strong and ordinary honest definitions, which characterize distality and NIP, respectively, in model theory \cite[Theorems~11 and~21]{CS}. This analogy led me to look for a counterexample to proper compression among NIP structures without distal expansions, focusing on the family of vertex neighborhoods in a disjoint union of expander graphs \cite{CSt,JNOdMS}. I thank Atticus Stonestrom and Riley Thornton for discussions related to this approach in spring and fall 2024, respectively.
I want to thank Michael Benedikt for introducing me to the restricted quantifier collapse (RQC) property \cite{BST,FZ,BH} in spring 2026, for unrelated reasons. This property holds in all distal theories and some stable theories. I realized that it suffices for proper compression and found a proof by hand for formulas without the finite cover property (see also Chase and Freitag \cite{CF} on equivalence queries). In summer 2026, using GPT-5.6, I obtained proper compression schemes of total size $O((\ell+1)2^\ell)$ for classes of Littlestone dimension $\ell$ (equivalently, for stable formulas \cite{CFModel}). Further direct prompting with GPT-5.6 and GPT-6 yielded little.
The distal and stable cases, at opposite extremes of NIP, suggested applying Simon's type decomposition theorem \cite{SimonDecomposition}. I was able to use it with GPT-6 to prove proper compression for formulas in NIP theories. As in \cite{CS,LS}, the proof required finite VC dimension for auxiliary families obtained by Boolean operations and projection. Using GPT-6 and then GPT-6.2, I was eventually able to remove this assumption (and to improve the proper compression bound to $O(\ell^2)$ for classes of Littlestone dimension $\ell$---I do not discuss this refinement here). I then substantially reworked the proof.

The author was partially supported by the NSF Research Grants DMS-2246598, DMS-2554164 and by the Alexander von Humboldt Foundation.

\end{document}